\documentclass[11pt]{article}
\usepackage{graphicx}
\usepackage{amssymb}

\begin{document}
\title{\textbf{Automorphisms of Certain Projective Bundles over Toric Varieties}} 
\author{Amassa Fauntleroy \\
Department of Mathematics\\
North Carolina State University\\
Raleigh, North Carolina 27695\\}
\date{}

\maketitle
\begin{abstract}
\noindent  The purpose of this note is to exhibit the automorphism group of a projective bundle P(E) over a simplicial toric variety $X$ when the bundle E is a direct sum of equvariant line bundles. This case is important in the study of moduli of complete intersections on toric varieties including projective spaces.  The main result is that the automorphism group of $P(E)$ is, up to a finite group, the semi-direct product of the automorphism group of the base and a certain subgroup of fiber preserving automorphisms. This structure is similar to the structure of the automorphism groups of rational surfaces $F_{n}$ (see \cite{f7} p. 425) Applications to moduli space constructions are indicated in special cases, including Del Pezzo surfaces and certain Calabi Yau m-folds.
\end{abstract}

\section{Introduction}

The automorphism group of a complete simplicial toric variety has been determined by several
authors. Among these are Demazure (\cite {Demazure}) in the special case where the variety is smooth and Cox in \cite{Cox} in
the general setting. The purpose of this investigation is to study the automorphism group of a projective bundle P(E) over a simplicial toric variety $X$ when the bundle E is a direct sum of equivariant line bundles. In this case P(E) is itself a toric variety. Though the groups involved here are not reductive, they are similar to them in the sense that they are generated by root subgroups. The systematic study of the possible root system structures has been made by Demazure in \cite {Demazure}. Even though the rich decomposition theory for the semisimple groups is not available, these roots are an efficient way to describe the actions of these groups on the algebraic structures involved when constructing moduli spaces. The main result of this note, \textbf{Theorem 1}, is that the automorphism group of $P(E)$ is, up to a well determined finite group, the semi-direct product of the automorphism group of the base and a certain subgroup of fiber preserving automorphisms. This structure is similar to the structure of the automorphism groups of rational
surfaces (see\cite {f7}).

We follow the notations and conventions of Baytrev and Cox (\cite {BC}]) and Oda (\cite {Oda}) for toric varieties. Let M be a free abelian group of rank d and N=Hom(M,$\mathbb{Z}$) the dual group. We denote by $M_\mathbb{R}$ and $N_\mathbb{R}$ the $\mathbb{R}$-scalar extensions of M and N. We denote by $\sigma$ a rational polyhedral cone in $N_\mathbb{R}$ and by $\Sigma$ a rational simplicial complete d-dimensional fan. (See Cox and Batryev \cite {BC} section 1 for details). Let $S = \mathbb{C}[x_1, ... , x_l]$ denote the homogeneous coordinate ring of $X(\Sigma)$ defined by Cox. The toric variety $X(\Sigma)$ defined by $\Sigma$ is the quotient of an open subset of the affine space Spec(S) by a torus, $D(\Sigma)$, whose action is defined by the combinatorial data of $\Sigma$ (see \cite{Cox} Theorem 2.1 ). The ring $S$ is graded in a natural way by the divisor class group $Cl(\Sigma) = A_{d-1}(X)$ of $X$. We consider $S$ with this fixed grading. We call  $\alpha$ the degree of an element $f(x_1, ... , x_l)$ in the graded subspace $S_{\alpha}$  of S and call the usual degree determined by integer exponents the polynomial degree of $f(x_1, ... , x_l)$. The group of graded automorphisms of $S$ will be denoted by $G$. This group maps onto the connected component of the automorphism group of $X$ with finite kernel (an isogeny). We will describe the automorphisms of toric varieties using this special presentation. An explicit description of $D(\Sigma)$ and of the finite group $Aut(X) / G $ can be found in Cox's paper \cite {Cox}.

Suppose $X$ is a toric variety and $L_1, ..., L_r$ are equivariant line bundles given by support functions
$h_{i}, i=1,...,r$ on $\Sigma$ (see\cite {Oda} section 2.1). Let $E$ denote the direct sum of the corresponding line bundles. The projective bundle $P(E)$ is a toric variety defined by the following simplicial fan:

Introduce the $\mathbb{Z}$-module $N'$ with a $\mathbb{Z}$ basis $n_2 ,\cdots,n_r$. Let $\tilde{N}:=N+N'$ and put $n_{1}= - n_{2}-\cdots- n_{r}$. Denote by  $\widetilde{\sigma}$  the image of each cone  $\sigma$  $\in$  $\Sigma$ under the $ \mathbb{R}$-linear map from $N_\mathbb{R}$ to $ \widetilde{N}_\mathbb{R}$ which sends y in $N_\mathbb{R}$ to $y - \sum_{j=1}^r h_j(e_i)n_j$. On the other hand let $\sigma_i'$ be the cone generated by $n_1,...n_{i-1},n_{i+1},...,n_r$ in $N'$ and let $\Sigma^{''}$ be the fan in $N'_\mathbb{R}$ consisting of the faces of  $\sigma_1',...,\sigma_s'$. Then the projective bundle $P(E)$ corresponds to the fan $\widetilde{\Sigma} :=\{ \tilde {\sigma} +\sigma' : \sigma \in \Sigma, \sigma' \in \Sigma^{''}\}$.  The homogeneous coordinate ring of $P(E)$ is then $\mathbb{C}[x_1,\cdots,x_l,y_1,\cdots,y_r]$

Next we recall the contravariant defnition of toric varieties $P_\Delta$ associated
to a lattice polyhedron $\Delta$. Start with $M$ as above and let $\Delta$ denote a convex d-dimensional cone
in $M_R$ with integral vertices(i.e., all vertices are in $M$ ). Let $S_\Delta$ denote the subring of
$C[x_0,x_1,...,x_d,x_1^{-1},...,x_d^{-1}]$ whose basis consist of monomials
$x_0^{m_0} \cdot x_1^{m_1} \cdots x_n^{m_n}$ where $m_0$ is positive and ($m_1/m_0$,...,$m_d/m_0$) is in $\Delta$. Then $P_\Delta$ is Proj($S_\Delta$). (See \cite {BC} for details and the connections between the two notions)

\section{Automorphism Groups of P(E)}

 Following Cox \cite{Cox} we use the homogeneous coordinate ring $S$ to describe the automorphism group of X. The crucial fact for us is that the full automorphism group is the isogeneous image of an affine algebraic
group whose connected component is the group of graded automorphisms of S. This group is particularly easy to
describe and contains the root subgroups mentioned in the introduction.

	 If $\alpha \in Cl(\Sigma)$ is the class of a polynomial generator $S$ and $S_\alpha$ is the corresponding eigenspace for the action of the torus $D(\Sigma)$, then a root through $\alpha$ consist of a map $\gamma$ from $G_a$ to G which acts as follows: there is a variable $x$ of weight $\alpha$ and a monomial m of weight $\alpha$ which is not divisible by $x$ such that $\gamma(t)$ is the indentity on all variables distinct from $x$ but sends that variable to $x + tm$.  Let us call a root reductive if it corresponds to a pair $(x , m)$ where the polynomial degree of m is 1. Otherwise, we say the root is unipotent. If the pair $(x , m)$ defines a reductive root so does the pair $(m,x)$. Such a pair gives rise to a homomorphism from $SL(2,\mathbb{C})$ to $G$. In fact, the only reductive subgroups which appear in the decomposition of $G / Rad_u (G)$ are of type A, where  $ Rad_u (G)$ denotes the unipotent radical of $G$.

	 Recall ( \cite {Cox} ) that $P_{\Delta}$ can be realized as the quotient of $U$=Spec$S$-V($\Sigma)$ by the action of the torus D whose character group is Cl($\Sigma$). Then if J is the normalizer of D in the automorphism group of $U$, G = Aut($S_g$) is the connected component of the identity and Aut($X$) = $J$/$D$. The group $G$ is the semidirect product of a Levi factor $G_l$ = $\bigoplus_{j=1}^sGL(r_j,\mathbb{C})$ and its unipotent radical $H$ generated by the unipotent roots above. The $r_{j}$ are the dimensions of the subspaces of the eigenspaces $S_{\alpha_{j}}$ spanned by the variables of degree $\alpha_{j}$ .

 Now let $X$ be a smooth toric variety and let $L_{1}, ... ,L_{r}$ be equivariant line bundles on $X$ defined by support functions $h_{j}, j=1,...,r$ on $\Sigma$ (see \cite{Oda}). Let $l = |\Sigma(1)|$ and let $e_{i}, i=1,...,l$ denote the vertices. Put $\widetilde{e_{i}}=e_{i}-\sum_{j=1}^l  h_{j} (e_{i}) n_{j}$ for $1\leqslant  i  \leqslant l$ where $n_{1}+ ...+n_{r}=0$,  $M=N^\ast$ and $N\otimes\mathbb{R}$ has basis among  $e_{i}, i=1,...,l$

The exact sequence $$0 \rightarrow M \stackrel {i} {\rightarrow} \mathbb {Z}^{l} \stackrel {\varphi_{0} } {\rightarrow} Pic(X) \cong A_{d-1}(X) \rightarrow 0 $$
defines $A_{d-1} (X)$. Recall $i(m)=(m(e_{1}), ..., m(e_{l}))$.  Let $M^{'}$=$M \oplus [\coprod_{i=2}^{r} \mathbb{Z}{n_{i}}]^{^\ast}$

Now define a sequence $$0 \rightarrow M^{'}\stackrel {i^{'}} {\rightarrow} \mathbb {Z}^{l} \oplus \mathbb {Z}^{r} \stackrel {\varphi} {\rightarrow} A_{d-1}(X) \oplus \mathbb{Z} \rightarrow 0 $$
where $$i{'} (m^{'} )=(m^{'}(\tilde{e}_{1}), ... ,m^{'}(\tilde{e}_l), m^{'}(n_{1}),...,m^{'}(n_{r}))$$ and $$\varphi (a_{1}, ...,a_{n},b_{1},...,b_{r})=(\sum_{i=1}^na_{i}\varphi_{0} (e_{i})-\sum_{j=1}^r b_{j}\alpha_{j}  ,  b_{1}+ ...+b_{r})$$
Here $\alpha_{j}$ is the class of $L_{j} \in A_{d-1}(X)  , 1 \leqslant j \leqslant r$.

\textbf{Lemma 2.1} There is a natural isomorphism $Pic(P(E)) \cong A_{d-1}(X) \oplus \mathbb {Z }$ where $E=L_{1} \oplus... \oplus L_{r}$. Moreover, the sequence $$ 0 \rightarrow M^{'} \rightarrow \mathbb {Z}^{l} \bigoplus\mathbb{Z}^{r}\rightarrow Pic(P(E))\rightarrow0$$ is exact.

\textbf{Proof} The first assertion is proven in \cite{Mav}. Now $(\varphi \circ i^{'}) (m , 0)=\varphi_{0} (i(m),0)=0$ by definition of $\varphi_{0}$ and $i$. Also, $i(0,n_{j})=(- h_{1}(e_{1})+ h_{j}(e_{1}) , ... , - h_{1}(e_{l}) + h_{j}(e_{l}),-1,0,...,0,1,0,...,0)$ and hence $\varphi (i(0,n_{j}))=(\sum_{i=1}^n (h_{j}(e_{i})-h_{1}(e_{i})\beta_{i} +\alpha_{1}-\alpha_{j} ,0)$.
Since $\alpha_{1}=\Sigma h_{1}(e_{i})\beta_{i}$ and $\alpha_{j}=\Sigma h_{j}(e_{i})\beta_{i}$, we see that $\varphi (i(0,n_{j}))=0$. 
Now $\varphi((0,...,1,0;0,...,0)=(\beta_{i},0)$, for $i=1, ..,l$ and $\varphi((0,...,,0;0,...,0,1,0,...,0)=(-\alpha_{j},1)$. Since the $\beta_{i}, i=1,...,l$ generate $A_{d-1}(X)$, $\varphi$ is surjective.

\textbf{Corollary 2.2} Let $ \{ x_{i}, y_{i}\}$ be the variables which generate the Cox coordinate ring of $[P(E)]$. If $f(x_{1}, ...,x_{n})$ has weight $\-\Sigma h_{j}(e_{i})\beta_{i}$ then $f(x_{1}, ...,x_{n})y_{j}$ has weight (0, 1) in $A_{d-1}(X)$.

\textbf{Proof} Since $\alpha_{j}=\Sigma h_{j}(e_{i})\beta_{i}$, we see from the lemma that $F(x_{1}, ...,x_{n})y_{j}$ has weight $(\alpha_{j}, 0)+(-\alpha_{j}, 1)=(0,1)$.

We want to apply the above observation to study the group $Aut(P(E))$. We introduce a bit more notation: $\pi:P(E )\rightarrow X$ denotes the projection map,  $\alpha_{i}=[L_{i}] \in A_{d-1}(X)$ and $G=Aut^{0}(P(E))$ the connected component of $Aut(P(E))$. Let $H$ be the subgroup of G consisting of bundle maps, i.e.,  if $h\in H$ then the induced map on $P(E)$ satisfies $\pi \circ h= \pi$. Let $K=Aut^{0}(X)$ be the connected component of $Aut(X)$.

\textbf{Theorem 2.3} There exists a monomorphism $i:K \rightarrow G$ such that for each $k \in K$, $k\circ \pi = \pi \circ i(k)$. Moreover, the subgroup $H$ is normalized by $i(K)$ and via $i$, $G$ is the semi-direct product of $H$ and $K$.

\textbf{Proof} Let $R$ and $S$ denote the Cox coordinate rings of $X$ and $P(E)$ respectively. Then $S$ is an $R$ algebra. Let $x_{1}, ... , x_{l}$ generate $R$ over $\mathbb{C}$ and let $y_{1}, ..., y_{r}$ generate $S$ over $R$. For $\sigma \in K$ extend the action of $\sigma$ to an automorphism of $S$  by defining $i( \sigma)(y_{j})=y_{j}$, $1\leqslant j \leqslant r$. Then $i:K \rightarrow G$ is one to one and the induced map on $P(E)$ satisfies $\pi \circ i(\sigma) = \sigma \circ \pi$.

Now suppose $h \in H$ and $p \in P(E)$. Then $$(\pi \circ i(\sigma))(h(p))=\sigma (\pi(h(p)) = \sigma  \pi(p)$$ so $\pi \circ i\sigma \circ h = \sigma \circ \pi$.

Note that $$ \pi( (i \sigma \circ h \circ i \sigma^{-1})(p))= \sigma(\pi (h \circ (i \sigma)^{-1})(p)=\sigma( \pi ((i\sigma)^{-1}(p))= \sigma^{-1} \circ \sigma (\pi(p))= \pi(p)$$ for all $p \in P(E)$. Thus $i \sigma \circ h \circ i \sigma^{-1} \in H$ and $i(K)$ normalize $H$.

If $\omega \in H \bigcap i(K)$, then $\pi (\omega(p))=\omega \circ \pi (p) =\pi(p)$. If $\omega=i(\sigma)$ with $\sigma \in K$, then $\omega \circ \pi (p) = k \circ \pi (p) = \pi(p)$ for all $p \in P(E)$. But $\pi$ is surjective so this means $\sigma$ and hence $\omega$ is the identity map.

Finally, to see that $K$ and $H$ generate all of $G$, it suffices to show that every root subgroup of $G$ lies in either $H$ or $K$. Now the lemma implies that every variable (i.e., polynomial generator) of $S$ is either in $R$ or one of the $y_{j}$ or in the $ \mathbb{C}$ - linear span of the $y_{j}$. If $(x,x^{D})$ corresponds to a root subgroup through the variable $x$, then the lemma implies the $x^{D} \in R$ (respectively the $\mathbb{C}$ - span of the $y_{j})$ when $x \in R$ (respectively the $\mathbb{C}$ - span of the $y_{j})$. Thus every root is either in $K$ or in $H$ as desired.

\textbf{Remarks}\begin{enumerate}
  \item  Neither $H$ nor $K$ is reductive or purely unipotent in general.
  \item  $Rad_{u}(G)$ is the semi-direct product of $Rad_{u}(H)$ and $Rad_{u}(K)$
\end{enumerate}

\section{Cayley Forms and Moduli  Spaces}

The notion of a Cayley form has been around for some time and is, in spirit,  a kind of Chow form.  We use a description given by A. Mavlyutov in \cite{Mav}. Let $f_{1}, \cdots ,f_{r}$ be homogeneous polynomials in $R$ with deg$f_{j}$ =$ [L_{j}]$. We assume from now on that the $L_{j}$ are ample. Then the Cayley form of $V(f_{1}, \cdots ,f_{r})= \{x \in X | f_{j}(x)=0,  1 \leqslant  j  \leqslant  r\}$ is $F=\sum f_{j}y_{j}$ in $S$. Moreover, the degree of $F$ is the class of the line bundle $O_{P(E)}(1)$ on P(E).

Recall that a subvariety of a toric variety P is quasi-smooth if its pre-image in the affine cone, Spec A ( A being the Cox coordinate ring), over P is smooth off the locus Z($\Sigma$) where the torus action fails to be proper (see \cite{Mav}). A complete intersection $V(f_{1}, \cdots ,f_{r})$ in a toric variety X is quasi-smooth if and only if the Cayley hypersurface $V(F)$  in P(E) defined by  the Cayley form F  is quasi-smooth.  

Now fix a smooth toric variety X, line bundles $L_{j}, 1 \leqslant j \leqslant r$ and sections $f_{j}$ of $L_{j}$.
Assume that $V(f_{1}, \cdots ,f_{r})$ is quasi-smooth. Put E = $L_{1} \bigoplus \cdots \bigoplus L_{r}$ and let F=$\Sigma f_{j}y_{j}$ be the corresponding Cayley form. Let $G = Aut^{0}P(E)$ and denote $i(K)$ (see \textbf{Theorem 2.3})  by $G_{b}$

\textbf{Lemma 3.1} Let $\omega$ be in $G$ with $\omega=\sigma\tau$, $\sigma \in G_{b}, \tau \in H$.  Then $\omega(F)$ is a Cayley form. Suppose $\omega (F) = \Sigma g_{j}y_{j}$. Let the ideals $I=\{f_{1} \cdots f_{r}\}$ and $J=\{g_{1} \cdots g_{r}\}$ be the coefficient ideals of $F$ and $\omega(F)$. Then $\sigma(I)=J$.

\textbf{Proof:} Note that if $\tau=1$, then $\omega(F)=\Sigma  \sigma (f_{j})y_{j}$ and the result is clear. It suffices then to establish the result for $\tau \in H$ and for this it is sufficient to examine root subgroups in $H$. 

If $\tau$ is a root subgroup in H then either\begin{itemize}
  \item there are indices t and s with $\tau (y_{t}) =y_{t} +h(x)y_{s}$ and deg($y_{t}$)=deg(h) + deg($y_{s}$) and $\tau(y_{j})=y_{j}$ for $j\neq t$ or
  \item there are indices t and s with $\tau (y_{t}) =y_{t} +k y_{s},  k \in\mathbb{C}$ and deg($y_{t}$) = deg($y_{s}$) and $\tau(y_{j})=y_{j}$ for $j\neq t$
\end{itemize}
For such $\tau$, $$\tau(F) = \sum_{j \neq t,s}f_{j}y_{j} + f_{t}(y_{t}+h(x)y_{s}) +f_{s}y_{s}= \sum_{j\neq s}f_{j}y_{j}  + (f_{s} +f_{t}h)y_{s}$$ Then the coefficient ideals I and J=$\tau$(I) are identical. Since the root subgroups generate H, the lemma follows.

Using the lemma we now give the connection to moduli constructions.

\textbf{Theorem 3.2:} Let $V_{1}$ and $V_{2}$ be smooth complete intersections in the toric variety X. Suppose each is of type $\beta$ in the Chow group of X and let  $V_{1}(E)$ and $V_{2}(E)$ be the corresponding Cayley hypersurfaces in P(E). Then the following are equivalent:

    1. There exists $\sigma \in Aut(X)$ with $\sigma(V_{1})=V_{2}$.

    2. There exists $\omega \in Aut(P(E))$ with $\omega(V_{1}(E))=V_{2}(E)$.

\textbf{Proof:} Assume that assumption 1 holds and let $f_{1}y_{1}+\cdots+f_{k}y_{k}$ be the Cayley form for $V_{1}$. Let $\hat{\sigma}$ denote the extension of $\sigma$ to P(E). Then $\hat{\sigma}(f_{1}y_{1}+\cdots+f_{k}y_{k})=\sigma(f_{1})y_{1}+\cdots+\sigma(f_{k})y_{k}$ is the Cayley form for $\hat{\sigma}(V_{2}(E))$ and 2 holds.

	Conversely, let $\omega (V_{1}(E))=V_{2}(E)$ and write $\omega = \hat{\sigma} \tau$ with $\hat{\sigma} $ the lift of $\sigma \in G_{b}$ and $\tau \in H$. The lemma implies $\sigma (I_{1}) = I_{2}$. Since $V_{1}$ is a smooth complete intersection, $I_{2}$ is generated by $\sigma (f_{1}) \cdots \sigma(f_{k})$ and so $\sigma$ maps $V_{1}$ to $V_{2}$. Q.E.D.
	
	In order to use the above result to construct moduli spaces we must restrict the class of complete intersections under consideration. Let $ \mathbf{\beta} = (\beta_{1}, \cdots ,\beta_{r})$, where $\beta_{j}  \in A_{d-1}(X)$, be a sequence of degrees. The sequence $\mathbf{\alpha}$ = $(deg ( x_{1}), \cdots ,deg ( x_{l}))$ will be called the fundamental sequence of degrees. The support $|\beta|$ of $\mathbf{\beta}$ is the set of distinct $\beta_{j}$ for some $j=1,...,r$. A complete intersection V is of type $\mathbf{\beta}$ if the defining functions $f_{j} ( x_{1}, \cdots, x_{l}) \in S_{\beta_{j}}$. We then write $V(\mathbf{\beta})$ for this subvariety.

	 Let us call a complete intersection $V(\mathbf{\beta})$ in $X(\Sigma)$ \textit{ nondegenerate}      if $|\mathbf{\beta}|\cap|\mathbf{\alpha}|=\emptyset$. Geometrically, this means (similar to the case of projective space) that V is not contained in any fundamental subspace ( 'hyperplane') of $X(\Sigma)$.  Next define the homogeneous coordinate ring of V to be the $\mathbb{C}$-algebra R=S/I where I is the ideal generated by the $f_{j}$. V  is said to be \textit{projectively normal} when R is an integrally closed domain. Let us continue to assume that $X(\Sigma)$ is smooth and let $L_{j} = O_{X}(\beta_{j})$  be the line bundle corresponding to $\beta{j} \in A_{d-1}(X) =Pic(X)$. As above we put $E = L_{1} \oplus \cdots \oplus L_{r}$. Let $V(\mathbf{\beta})$ be a complete intersection of type $\mathbf{\beta}$	and denote by $V(\mathbf{\beta},E)$ the hypersurface in P(E) defined by the Cayley form of $V(\mathbf{\beta})$. 
	
	\textbf{Lemma 3.3} If $V(\mathbf{\beta})$ is smooth and nondegenerate then the Cayley hypersurface  $V(\mathbf{\beta},E)$ is smooth, nondegenerate and projectively normal in P(E).

	\textbf{Proof}: Smoothness is clear (see also \cite{Mav} ). Since $V(\mathbf{\beta})$ is smooth, the defining polynomials $f_{1},\cdots, f_{r}$ have no common factor and so the Cayley form is irreducible. It follows from standard commutative algebra arguments (see eg \cite{f8}) that $V(\mathbf{\beta},E)$ is projectively normal.
	
	The arguments given in \cite{f8} can now be used to prove 
	
	\textbf{Theorem 3.4:} The moduli space of nondegenerate smooth complete intersections $V(\mathbf{\beta})$ of type $\mathbf{\beta}$ on the smooth toric variety  $X(\Sigma)$ is isomorphic to the quotient space of the open subset of properly stable points of  $H^{0}(P(E), O_{P(E)}(1))$ by the action of the group $G=Aut_{g}(S)$, where S is the homogeneous coordinate ring of P(E).
	
	\textbf{Proof:} Let $V_{1}(\mathbf{\beta},E)$ and $V_{2}(\mathbf{\beta},E)$ be the Cayley hypersurfaces of two isomorphic smooth complete intersections of type $\mathbf{\beta}$ in X .  By projective normality, there is an induced isomorphism of their homogeneous coordinate rings which, by nondegeneracy, lifts to an automorphism of S. It is straight forward to check that the element of G so defined maps $V_{1}(\mathbf{\beta},E$ to $V_{2}(\mathbf{\beta},E)$ (see \cite {f7} Theorem 4.3). The converse is just Theorem 3.2 above.
	
	\textbf{Remark}: If $\pi : V(\mathbf{\beta},E) \rightarrow X$ is the restriction of the projection map from P(E) to X, then $V$ is precisely the locus in X where  $\pi$ fails to have finite fibers-the 'exceptional locus' of $\pi$.
	
	\textbf{Examples:}
	
	1.  Fix a positive integer m and let V be the subvariety of $P^{m+r}$ defined by the vanishing of the r homogeneous forms $f_{j}$ of degrees $d_{j}$ for $ 1 \leqslant j\leqslant r$. If the resulting V is a complete intersection of dimension m and the sum of the degrees is m+r+1, then V is a Calabi- Yau m-fold. When m=3, r=2 and $d_{1}=2, d_{2}=4$ we get in particular a Calabi-Yau 3-fold. Here $E= O_{P^{5}}(2)\oplus O_{P^{5}}(4)$. 
	
	By the above we see that the space of Cayley forms is 
	 $$H^{0}(P(E), O_{P(E)}(1)) = \{f \cdot y_{1}+g \cdot y_{2} : f \in H^{0}(P^{5},O_{P^{5}}(2)),  g \in H^{0}(P^{5},O_{P^{5}}(4))\}$$
	The dimension of the space is 147. The automorphism group, G, of the homogeneous coordinate ring of P(E) is the isogeneous image of the semi-direct product $$[GL(6,\mathbb{,C})\times \mathbb{C}^{*} \times \mathbb{C}^{*}]\cdot H$$ where the unipotent radical H is isomophic as vector group to $H^{0}(P^{5}, O_{P^{5}}(2))$. Cayley forms F and G corresponding to $(f,g)$ and $(f^{'},g^{'})$ are isomorphic if and only if they are in the same G orbit.
	The group $[GL(6,\mathbb{,C})\times \mathbb{C}^{*} \times \mathbb{C}^{*}]\cdot H$  has dimension 59. Taking account of the torus D($\Sigma)$= $\mathbb{C}^{*} \times \mathbb{C}^{*}$ which is acting to produce the quotient $X(\Sigma)$=P(E), the correct dimension to use is 57. Thus the moduli space has dimension 90. If we use the more classical approach, the forms of degrees 2 and 4 on    $P^{5}$ have (projective) dimensions 20 and 125 respectively. The dimension of $PGL(6,\mathbb{C})$ is 35 and again the dimension of moduli is 90.
	
	2. Let P(1,1,2,3) be the weighted projective space and S=$\mathbb{C}[x,y,z,w]$ the homogeneous coordinate ring where deg(x)=deg(y)=1, deg(z)=2 and deg(w)=3. The automorphism group,  $G$, in this case is isogenous to $[GL(2,\mathbb{,C})\times \mathbb{C}^{*} \times \mathbb{C}^{*}]\cdot H$ where the unipotent radical, H, is generated by two types of roots. The first type sends z to z+h(x,y) where h is homogeneous of degree 2 and leaves x, y and w fixed. The second type maps w to w+g(x,y,z) where g is homogeneous of  degree 3 in S and leaves x, y and z fixed. 
	Hypersurfaces of degree 6 in S define Del Pezzo surfaces. Using the action of the unipotent radical of $G$ one can assume such a form  is given by the equation $$ aw^{2} + bz^{3} +zf_{4}(x,y) +f_{6}(x,y)=0$$ where$f_{4}(x,y)$ and $f_{6}(x,y)$ have degree 4 and 6 respectively. Next, using the subgroup corresponding to
$\mathbb{C}^{*} \times \mathbb{C}^{*}$ one can assume a=b=1. The moduli now depend on the pair of binary forms $f_{4}(x,y)$ and $f_{6}(x,y)$. The group $GL(2,\mathbb{,C})$ acts on these pairs and the moduli problem is now a classical one (see \cite{mum} ). In fact, using models of Del Pezzo fibrations constructed by Grinenko \cite{grin}, it seems likely that these group actions also play an important role in the moduli of these 3-folds.	

In this case, H is itself a semi-direct product, not a vector group. it is possible in general to describe the structure of the unipotent radical H using a normal series of subgroups determined by suitable orderings of the roots. Note for example, that in this case, the subgroup generated by roots of type 2 is normalized by the subgroup generated by roots of the first type and each type generates a vector group. The dimension of H is 9 in this case. A minimal toric desingularization (in the sense of the number of vertices that needs to be added) will yield an automorphism group whose unipotent radical has dimension 11. It would be interesting to know what algebro-geometric information is being carried by the orbits of these extra 2 dimensions.

\end{document}